\documentclass{article}
\usepackage{amsmath,amsthm,amsfonts, latexsym}

\newtheorem{thm}{Theorem}[section]

\newtheorem{lem}[thm]{Lemma}
\newtheorem{prop}[thm]{Proposition}
\newtheorem{cor}[thm]{Corollary}

\theoremstyle{definition}
\newtheorem*{defin}{Definition}
\newtheorem{ex}{Example}

\begin{document}

\begin{center}
\Large \textbf{Semilocal Generic Formal Fibers} \\
\vspace{0.2in}
\large P. Charters and S. Loepp

\end{center}

\begin{abstract}
Let $T$ be a complete local ring and $C$ a finite set of incomparable prime
ideals of $T$.  We find necessary and sufficient conditions
for $T$ to be the completion of an integral domain whose
generic formal fiber is semilocal with maximal ideals the
elements of $C$.  In addition, if char$T = 0$, we give necessary
and sufficient conditions for $T$ to be the completion of
an excellent integral domain whose generic formal fiber is
semilocal with maximal ideals the elements of $C$.
\end{abstract}

\section{\textbf{\sc{Introduction}}}

If $A$ is a local integral doman with maximal ideal $M$, quotient
field $K$, and $M$-adic completion $\hat{A}$, then $\mbox{Spec}
(\hat{A} \otimes_{A} K)$ is called the generic formal fiber of $A$.
Note that there is a one-to-one correspondence between the
elements of the generic formal fiber of $A$ and the inverse
image of the ideal $(0)$ under the map $\mbox{Spec}\hat{A}
\rightarrow \mbox{Spec}A$.  In light of this correspondence,
if $Q \in \mbox{Spec} \hat{A}$ and $Q \cap A = (0)$, we
will say that $Q$ is in the generic formal fiber of $A$.
Furthermore,
if the ring $\hat{A} \otimes_{A} K$ is semilocal with
maximal ideals $P_{1} \otimes_{A} K,P_{2} \otimes_{A} K, \ldots ,
P_{n} \otimes_{A} K$, then we will say that the generic formal
fiber of $A$ is semilocal with maximal ideals $P_{1},P_{2},
\ldots ,P_{n}$.

Because the standard integral domains we study have generic formal fibers that
are far from semilocal, at first glance one might
guess that noncomplete domains
posessing a semilocal generic formal fiber do not exist.
However, in \cite{loepp}, it was shown that such rings
do exist
and perhaps even more surprisingly,
in \cite{loepp2}, it was shown that these integral
domains can be constructed to be excellent.  In this
paper, we show that these domains are more
plentiful than one might suspect (both
in the nonexcellent and excellent case).

In section \ref{nonexcellent}, we
characterize which complete local rings are
completions of integral domains posessing a semilocal
generic formal fiber.  Specifically, suppose
 $(T,M)$ is a complete local ring, and $G \subseteq \mbox{Spec}T$
such that $G$ is nonempty and
the number of maximal elements of $G$ is finite.
We show that there exists a local domain $A$ such
that $\hat A = T$ and the generic formal fiber of $A$ is exactly $G$
if and only if $T$ is a field and $G = \{(0)\}$
or the following conditions hold.
\begin{enumerate}
\item $M \notin G$, and $G$ contains all the associated primes of $T$
\item If $Q \in G$ and $P\in \mbox{Spec}T$ with $P \subseteq Q$ then $P \in G$
\item If $Q \in G$ then $Q \cap \mbox{ prime subring of }T = (0)$
\end{enumerate}
It is easily seen that the above three conditions are
necessary and so the bulk of the proof
is dedicated to showing that the conditions are sufficient.
It is worth pointing out that the three conditions in
our theorem
are relatively weak, and so in some sense most complete
local rings can be realized as the completion of an integral domain
whose generic formal fiber is semilocal where
the maximal ideals can be prescribed.

In section \ref{excellent} we tackle the analogous
version of the above problem where we require the
additional condition that $A$ be excellent.
We are successful in characterizing the
complete local rings of characteristic zero that are
completions of excellent integral domains posessing a semilocal
generic formal fiber.
Specifically, let $(T,M)$ be a complete local ring containing the integers.
Let $G
\subseteq \mbox{Spec}T$ such that $G$ is nonempty and the number of
maximal elements of $G$ is finite.  We show there exists an
excellent local domain $A$ with $\hat{A} = T$ and
such that $A$ has generic formal fiber exactly $G$ if
and only if $T$ is a field and $G = \{(0)\}$
or the following conditions hold.
\begin{enumerate}
\item $M \notin G$, and $G$ contains all the associated primes of $T$
\item If $Q \in G$ and $P\in \mbox{Spec}T$ with $P \subseteq Q$ then $P \in G$
\item If $Q \in G$, then $Q \cap \mbox{ prime subring of }T = (0)$
\item $T$ is equidimensional
\item $T_{P}$ is a regular local ring for all maximal elements $P \in
G$.
\end{enumerate}
Showing that the above five conditions are necessary, although
maybe not immediately obvious, is relatively short.  Our proof,
then, will focus on proving that they are sufficient.

For both theorems, to show that the respective conditions
are sufficient we construct the desired integral
domain $A$.  Our construction is based on
the on the techniques used in \cite{loex}.
We start with the prime subring of $T$, localized at the appropriate prime
ideal.  We then successively adjoin elements of $T$ to this ring in order
to get our final result.  Naturally, we must be careful which elements we
choose to adjoin.  For example, we must avoid the zero divisors of $T$,
so that $A$ will be an integral domain.  We must also
avoid nonzero elements of prime ideals that we wish
to be in the generic formal fiber of $A$.
We will adjoin enough elements of $T$ to our domain
$A$ so that if $I$ is a finitely generated ideal of $A$
then $IT \cap A = I$.  Furthermore, we will
be adjoining elements of $T$ until we have obtained the property that for every
ideal $J$ of $T$ such that $J \not\subseteq P$ for all $P \in G$,
our ring contains a nonzero element of
every coset in the ring $T/J$.  Thus our ring will satisfy the property
that if $J$ is an ideal of $T$ where $J \not\subseteq P$ for all $P
\in G$, then the map $A \rightarrow T/J$ is onto.  In
particular, this means that $A \rightarrow T/M^{2}$ is onto.
This fact, along with the condition that $IT \cap A = I$ for
every finitely generated ideal $I$ of $A$ will force
the completion of $A$ to be $T$.  Moreover, what is
also interesting about the condition that
$A \rightarrow T/J$ be onto
is that it turns out if $T$
contains the integers then it will force $A$ to be excellent.  By adjoining
nonzero elements of each ideal $J$ where
$J \not\subseteq P$ for all $P \in G$ while avoiding
nonzero elements of the prime ideals contained in $G$,
 we also ensure that the generic formal fiber of $A$ is
exactly $G$.

All rings in this paper are to be assumed commutative with
unity.  If we say a ring is local, we mean it is a Noetherian ring
with one maximal ideal.  The term quasi-local will be reserved
for a ring with one maximal ideal that need not be Noetherian.
We will use $c$ to denote the cardinality of the real numbers.

\section{\textbf{\sc{The Construction}}} \label{construction}



We now begin the construction of our integral domain $A$.
The following proposition is Proposition 1 from \cite{H2}.
It will be used to show that the ring $A$ we construct has the
desired completion.

\begin{prop}  \label{prop1}
If $(A, M \cap A)$ is a quasi-local subring
of a complete local ring $(T,M)$, the map $ A \longrightarrow T/M^2$
is onto and $IT \cap A = I$ for every finitely generated
ideal $I$ of $A$, then $A$ is Noetherian and the natural
homomorphism $\hat{A} \longrightarrow T$ is an isomorphism.
\end{prop}

Although Lemma \ref{card} is well-known, we will use it
repeatedly.  So, we state it here without proof.

\begin{lem}\label{card}
Let $T$ be an integral domain and $I$ a nonzero ideal
of $T$.  Then, $|I| = |T|$.
\end{lem}


\begin{lem} \label{tp}
Let $(T,M)$ be a complete local ring of dimension
at least one. Let $P$ be a nonmaximal prime ideal of $T$.
Then, $|T/P| = |T| \geq c$.
\end{lem}

\begin{proof}
Clearly, $T/P$ is reduced.  Furthermore, since $T$ is complete and dim$T
\geq 1$,
$T/P$ is complete and dim$(T/P) \geq 1$, as $P$ is nonmaximal.
Since $T/P$ is reduced, complete and dim$(T/P) \geq 1$, we have
$|T/P| \geq c$.  But clearly $|T/P| \leq |T|$, so $|T| \geq c$.
Now, define a map $f: T \rightarrow
\prod_{i = 1}^{\infty} T/M^{i}$ by
$f(t) = (t + M,t + M^{2}, t + M^{3}, \ldots)$.  It is easy
to see that $f$ is injective and so $|T| = \sup\{c,|T/M|\}$.
Now, $|T/P| \leq |T|$ and $|T/P| \geq \sup\{c,|T/M|\} = |T|$,
so $|T/P| = |T|$ as desired.
\end{proof}

Armed with the previous two lemmas, we can now prove the following
critical lemma.
It will be used to adjoin elements to a specific
subring of $T$ so that the resulting ring maintains certain
properties of the original subring.

\begin{lem} \label{notun}
Let $(T,M)$ be a complete local ring such that $\mbox{dim}T \geq 1$,
 $C$ a finite set of nonmaximal prime
ideals such that no ideal in $C$ is contained in another ideal of $C$,
and $D$ a subset of $T$ such that
$|D| < |T|$.  Let $I$ be an ideal of $T$ such that $I \not\subseteq
P$ for all $P \in C$.  Then $I \not\subseteq \bigcup \{ r+P | r \in D$, $P \in
C\}$.
\end{lem}

\begin{proof}
Let $C = \{P_{1}, P_{2}, \ldots, P_{n}\}$.  From the Prime Avoidance
Theroem, we know that $I \not\subseteq
\bigcup_{i=1}^{n} P_{i}$.  Let $x \in I$, $x \notin
\bigcup_{i=1}^{n}P_{i}$.  Define a family of maps
$f_{i}: P_{i} \times D \to T$ for every $P_{i} \in C$ as follows.  Let
$(P_{i},r)
\in P_{i} \times D$.  If $r + P_{i} \notin (x + P_{i})(T/P_{i})$, define
$f_{i}(P_{i},r) = 0$.  Otherwise, it must be the case that $r + P_{i} = (x +
P_{i})(s_{i} + P_{i})$ for some $s_{i} \in T$, so choose one such
$s_{i}$ and
define $f_{i}(P_{i},r) = s_{i}$.  Now let $S_{i} = \mbox{Image}f_{i}$.
Note that we then have the inequality $|S_{i}| \leq |D| < |T| = |T/P_{i}|$.

First, suppose $n = 1$.  Then $|S_{1}| \leq |D| < |T| = |T/P_{1}|$.
So, there exists $t \in T$ such that $t + P_{1} \neq s + P_{1}$ for
all $s \in S_{1}$.  Now, if $xt \in \bigcup \{r + P_{1} |  r \in D\}$,
then $xt + P_{1} = r + P_{1}$ for some $r \in D$.  But then
$r + P_{1} \in (x + P_{1})(T/P_{1})$, so $r + P_{1} =
(x + P_{1})(s + P_{1})$ for some $s \in S_{1}$.  So, we have
$(x + P_{1})(t + P_{1}) = r + P_{1} = (x + P_{1})(s + P_{1})$
which implies that $t + P_{1} = s + P_{1}$, a contradiction.
It follows that the lemma holds if $n = 1$.


If $n > 1$, we claim that
$|T/P_{i}| = \left| \frac{P_{i} + \bigcap_{j = 1, j \neq i}^{n}
P_{j}}{P_{i}}\right|$.
Notice that since $T/P_{i}$ is an integral
domain, this is true by Lemma \ref{card}
if we can simply show that $\frac{P_{i} + \bigcap_{j = 1, j \neq i}^{n}
P_{j}}{P_{i}}$ is not the zero ideal of $T/P_{i}$.  Suppose
that this were not true.  Then it must be the case that $\cap_{j =1,
j \neq i}^{n}P_{j} \subseteq P_{i}$.  We know, however, that since no $P_{i}$
is contained in any other ideal in $C$ this cannot happen. Hence
$\frac{P_{i} + \bigcap_{j = 1, j \neq i}^{n}
P_{j}}{P_{i}}$ is not the zero ideal of $T/P_{i}$, and it follows
that $|D| < \left| \frac{P_{i} + \bigcap_{j = 1, j \neq i}^{n}
P_{j}}{P_{i}}\right|$.  Thus there exists a $t_{i} \in \cap_{j = 1,
j\neq i}^{n}P_{j}$ such that $t_{i} + P_{i} \neq s_{i} + P_{i}$ for
all $s_{i} \in S_{i}$ and for all $i = 1,\ldots,n$.
We claim that
$x\sum_{j=1}^{n}t_{j} \notin \bigcup \{ r + P_{i} | r \in D$, $P_{i} \in
C\}$.
To see this, suppose that $x\sum_{j=1}^{n}t_{j} \in \bigcup
\{ r + P_{i} | r \in D$, $P_{i} \in C\}$.  Then $x\sum_{j = 1}^{n}t_{j} + P_{i}
 = r + P_{i}$ for some $P_{i} \in C$, $r \in D$.  But this means that
$xt_{i} + P_{i} = r + P_{i}$, implying that $(x + P_{i})(t_{i} +
P_{i}) = r + P_{i}$ and thus $r + P_{i} \in (x + P_{i})(T/P_{i})$.
But then $(x + P_{i})(t_{i} + P_{i}) = r + P_{i} = (x +
P_{i})(s_{i} + P_{i})$ for some $s_{i} \in S_{i}$.  Thus $t_{i} +
P_{i} = s_{i} + P_{i}$ for some $s_{i} \in S_{i}$, a contradiction.
\end{proof}

\begin{defin}
Let $(T,M)$ be a complete local ring, and $C$ a set of prime
ideals of $T$.  Suppose that $(R, R \cap M)$ is a quasi-local subring
of $T$ such that $|R| < |T|$ and $R \cap P = (0)$ for every $P \in
C$.  Then we call $R$ a small $C$-avoiding subring of $T$ and will
denote it by $SCA$-subring.
\end{defin}

$SCA$-subrings will be essential in our proof.  If
$R$ is an $SCA$-subring of $T$ then note that if we choose our set $C$ such
that
the associated primes of $T$ are contained in
prime ideals in $C$, then the
condition $R \cap P = (0)$ for all $P \in C$
implies that $R \cap Q = (0)$ for every $Q \in
\mbox{Ass}T$, and thus $R$ contains no zero divisors of $T$ -
certainly a condition that any domain we might wish to
construct must enjoy.  Furthermore, this condition will ensure that
the prime ideals of $C$ are
in the generic formal fiber of our final domain $A$.
It is worth noting too that the condition $|R| < |T|$ implies that
$|R| < |T/P|$ for all nonmaximal prime ideals $P$ of $T$
from Lemma \ref{tp}.  This
cardinality condition will allow us to adjoin an
element to $R$ so that the resulting ring will
not contain zero divisors of $T$ or nonzero elements of
the prime ideals in $C$.

Recall that one property that we would like our constructed ring,
call it $A$, to possess, is that
if $J$ is an ideal of $T$ with $J \not\subseteq P$ for
all $P \in C$, then the map $A
\rightarrow
T/J$ is onto.  Lemma \ref{inmap} allows us to adjoin
an element of a coset of $T/J$, which eventually
will force our ring $A$ to satisfy this property.
The proof of Lemma \ref{inmap} closely parallels the proof of
Lemma 3 in \cite{LR} and Lemma 3 in \cite{loex}.

\begin{lem} \label{inmap}
Let $(T,M)$ be a complete local ring of dimension at least one.  Let
$C$ be a finite set of
nonmaximal prime
ideals of $T$ such that no ideal in $C$ is contained in any other
ideal in $C$.  Let
$J$ be an ideal of $T$ such that $J \not\subseteq P$ for all $P \in C$.
Let $R$ be an $SCA$-subring of $T$ and $u + J \in
T/J$.  Then there exists an infinite $SCA$-subring $S$ of $T$ such
that $R \subseteq S \subseteq T$ and $u + J$ is in the image of the
map $S \to T/J$.  Moreover, if $u \in J$, then $S \cap J \neq (0)$.
\end{lem}

\begin{proof}
Let $P \in C$.  Let $D_{(P)}$ be a full set of coset representatives
of the cosets $t + P$ that make $(u + t) + P$ algebraic over $R$.
Note that as $|R| < |T|$ and $|T| \geq c$, we have $|D_{(P)}| < |T|$.  Let
$D = \bigcup_{P \in C}D_{(P)}$, and note that $|D| < |T|$.  Now
use Lemma \ref{notun}
with $I = J$ to find an $x \in J$ such
that $x \notin \bigcup \{ r + P | r \in D$, $P \in C\}$.
We claim that $S = R[u+x]_{(R[u+x]\cap M)}$ is
the desired $SCA$-subring.  It is easy to see that $|S| < |T|$.
Now suppose that $f \in R[u+x] \cap P$ for some $P \in C$.  Then $f =
r_{n}(u+x)^{n} + \cdots + r_{1}(u+x) + r_{0} \in P$ where $r_{i} \in
R$.  But we chose $x$ such that $(u+x) + P$ is transcendental over
$R$.  Therefore $r_{i} \in R \cap P = (0)$ for every $i =
1,2,\ldots ,n$ and it follows that that $f = 0$.  So $S \cap P = (0)$
and we have $S$ is an $SCA$-subring.  Note further that if $u \in J$,
then $u+x \in J$.  Since $(u+x) + P$ is transcendental over $R$, it
must be the case that $u + x \neq 0$.  It follows that $S \cap J \neq
(0)$.
\end{proof}

The following lemma will help us ensure that $IT \cap A =
I$ for every finitely generated ideal of $A$.  Recall that this is a
necessary condition in order to be able to use Proposition \ref{prop1}.
The proof of Lemma \ref{scaex} resembles
that of Lemma 6 in \cite{loepp2}, as well as that of Lemma 4 in
\cite{loex}.

\begin{lem} \label{scaex}
Let $(T,M)$ be a complete local ring of dimension at least one.  Let
$C$ be a finite set of
nonmaximal prime ideals of $T$ such that if $Q \in \mbox{Ass}T$ then
$Q \subseteq P$ for some $P \in C$ and no ideal in $C$ is contained in
any other ideal in $C$,
and let $R$ be an $SCA$-subring of $T$.  Suppose that $I$ is a
finitely generated ideal of $R$ and $c \in IT \cap R$.  Then there
exists an $SCA$-subring $S$ of $T$ such that $R \subseteq S \subseteq
T$ and $c \in IS$.
\end{lem}

\begin{proof}
We will induct on the number of generators of $I$.  Suppose $I = aR$.
Now if $a = 0$, then $c = 0$ and thus $S = R$ is the desired
$SCA$-subring of $T$.  Thus consider the case where $a \neq 0$.  In
this case, $c = au$ for some $u \in T$.  We claim that $S = R[u]_{(R[u]
\cap M)}$ is the desired $SCA$-subring.  To see this, first note that
$|S| < |T|$.  Now let $f \in R[u] \cap P$ where $P \in
C$.  Then $f = r_{n}u^{n} + \cdots + r_{1}u + r_{0} \in P$.
Multiplying through by $a^{n}$, we get
$a^{n}f = r_{n}(au)^{n} + \cdots + r_{1}a^{n-1}(au) + r_{0}a^{n}$ and
it follows that
$a^{n}f = r_{n}c^{n} + \cdots + r_{1}a^{n-1}c +
r_{0}a^{n} \in P \cap R = (0)$.
Now, $a \in R$, $R \cap P = (0)$ for every $P \in C$, and
all associated prime ideals of $T$ are contained
in an element of $C$.  It follows that $a$ is not
a zero divisor in $T$.
It must be the case then that $f = 0$, giving us that $S$ is an $SCA$-subring
of $T$.  Thus we have proven the base case, when $I$ is principal.

Now let $I$ be an ideal of $R$ that is generated by $m > 1$ elements,
and suppose that the lemma holds true for all ideals of $R$ generated
by $m - 1$ elements.  Let $I = (y_{1}, \ldots , y_{m})R$.  Then $c =
y_{1}t_{1} + y_{2}t_{2} + \cdots + y_{m}t_{m}$ for some $t_{1},t_{2}
\ldots ,t_{n} \in
T$.  By adding $0$, note that we then have the
equality $c = y_{1}t_{1} + y_{1}y_{2}t - y_{1}y_{2}t + y_{2}t_{2} +
\cdots + y_{m}t_{m} = y_{1}(t_{1} + y_{2}t) + y_{2}(t_{2} - y_{1}t) +
y_{3}t_{3} + \ldots + y_{m}t_{m}$ for any $t \in T$.
Let $x_{1} = t_{1} + y_{2}t$ and
$x_{2} = t_{2} - y_{1}t$ where we will choose the element $t$ later.
Now, let $P \in C$.  If $(t_{1} + y_{2}t) + P = (t_{1} + y_{2}t') +
P$, then it must be the case that $y_{2}(t - t') \in P$.  But $y_{2}
\in R$, $R \cap P = (0)$ and $y_{2} \neq 0$, so we have $t - t' \in
P$.  Thus $t + P = t' + P$.  The contrapositive of this result
indicates that if $t + P \neq t' + P$, then $(t_{1} + y_{2}t) + P \neq
(t_{1} + y_{2}t') + P$.  Let $D_{(P)}$ be a full set of coset
representatives of the cosets $t + P$ that make $x_{1} + P$ algebraic
over $R$.  Let $D = \bigcup_{P \in C} D_{(P)}$. Note that $|D| < |T|$.
Now we can use Lemma \ref{notun}
with $I = T$ to find an element $t
\in T$ such that $x_{1} + P$ is transcendental over $R$ for every $P
\in C$.  It can be easily shown (as in the proof of
Lemma \ref{inmap}) that $R' = R[x_{1}]_{(R[x_{1}]\cap
M)}$ is an $SCA$-subring of $T$.  Now let $J = (y_{2}, \ldots ,
y_{m})R'$ and $c^{*} = c - y_{1}x_{1}$.  It is then the case that
$c^{*} \in JT \cap R'$, so we can use our induction assumption to draw
the conclusion that there exists an $SCA$-subring $S$ of $T$ such
that $R' \subseteq S \subseteq T$ and $c^{*} \in JS$.  Thus $c^{*} =
y_{2}s_{2} + \cdots + y_{m}s_{m}$ for some $s_{1},s_{2}, \ldots
,s_{n} \in S$.  It follows
that $c = y_{1}x_{1} + y_{2}s_{2} + \cdots + y_{m}s_{m} \in IS$, and
thus $S$ is the desired $SCA$-subring.
\end{proof}

\begin{defin}
Let $\Omega$ be a well-ordered set and $\alpha \in \Omega$.
We define $\gamma(\alpha) = \sup \{\beta \in \Omega |
\beta < \alpha\}$.
\end{defin}

Lemma \ref{sca} allows us to put many of our desired conditions
together.  We note here that the proof of Lemma \ref{sca} is based on
the proof of Lemma 12 in \cite{loepp2}.

\begin{lem} \label{sca}
Let $(T,M)$ be a complete local ring of dimension at least one.  Let
$J$ be an ideal of $T$ with $J \not\subseteq P$ for all $P \in C$,
where $C$ is a finite set of nonmaximal ideals of $T$ such that if $Q \in
\mbox{Ass}T$ then $Q \subseteq P$ for some $P \in C$ and no ideal
in $C$ is contained in any other ideal in $C$,
and let $u +
J \in T/J$.  Suppose $R$ is an $SCA$-subring.  Then there exists an
$SCA$-subring $S$ of $T$ such that
\begin{enumerate}
\item $R \subseteq S \subseteq T$
\item If $u \in J$, then $S \cap J \neq (0)$
\item $u + J$ is in the image of the map $S \to T/J$
\item For every finitely generated ideal $I$ of $S$, we have $IT \cap
S = I$.
\end{enumerate}
\end{lem}

\begin{proof}
We first apply Lemma \ref{inmap} to find an infinite
$SCA$-subring $R'$ of $T$ such that $R \subseteq R' \subseteq T$,
$u + J$ is in the image of the map $R' \to T/J$, and if $u \in J$
then $R' \cap J \neq (0)$.  We will construct the desired $S$ such
that $R' \subseteq S \subseteq T$ which will ensure that the first
three conditions of the lemma hold true.  Now let
\[
\Omega = \{ (I,c) | I \mbox{ is a finitely generated ideal of } R'
\mbox{ and } c \in IT \cap R'\}
\]
Letting $I = R'$, we can see that $|\Omega | \geq |R'|$.  But then
since $R'$ is infinite, the number of finitely generated ideals of
$R'$ is $|R'|$, and therefore $|R'| \geq |\Omega |$, giving us the
equality $|R'| = |\Omega |$.  Moreover, as $R'$ is an $SCA$-subring
of $T$, we have $|\Omega | = |R'| < |T|$.  Well order $\Omega$ so
that it does not have a maximal element and let $0$ denote its first
element.  We will now inductively define a family of $SCA$-subrings
of $T$, one for each element of $\Omega$.  Let $R_{0} = R'$, and
let $\alpha \in \Omega$.
Assume that $R_{\beta}$ has been defined for all
$\beta < \alpha$.  If $\gamma (\alpha) < \alpha$ and
$\gamma(\alpha) = (I, c)$, then define $R_{\alpha}$ to be the
$SCA$-subring obtained from Lemma \ref{scaex}.  In this manner,
$R_{\alpha}$ will have the properties that $R_{\gamma(\alpha)}
\subseteq R_{\alpha} \subseteq T$, and $c\in IR_{\alpha}$.  If
$\gamma(\alpha) = \alpha$, define $R_{\alpha} = \bigcup_{\beta <
\alpha}R_{\beta}$.  Note that in both cases, $R_{\alpha}$ is an
$SCA$-subring of $T$.  Now let $R_{1} = \bigcup_{\alpha \in
\Omega}R_{\alpha}$.  We know that $|\Omega | < |T|$ and
$|R_{\alpha}| < |T|$ for every $\alpha \in \Omega$, and thus
$|R_{1}| < |T|$ as well.  Moreover, as $R_{\alpha} \cap P = (0)$ for
every $P \in C$ and every $\alpha \in \Omega$, we have $R_{1} \cap
P = (0)$ for every $P \in C$.  It follows that $R_{1}$ is an
$SCA$-subring.  Furthermore, notice that if $I$ is a finitely
generated ideal of $R_{0}$ and $c \in IT \cap R_{0}$, then $(I,c) =
\gamma(\alpha)$ for some $\alpha \in \Omega$ with $\gamma(\alpha) <
\alpha$.  It follows from the construction that $c \in IR_{\alpha}
\subseteq IR_{1}$.  Thus $IT \cap R_{0} \subseteq IR_{1}$ for every
$I$ a finitely generated ideal of $R_{0}$.

Following this same pattern, build an $SCA$-subring $R_{2}$ of $T$
such that $R_{1} \subseteq R_{2} \subseteq T$ and $IT \cap R_{1}
\subseteq IR_{2}$ for every finitely generated ideal $I$ of $R_{1}$.
Continue to form a chain $R_{0} \subseteq R_{1} \subseteq R_{2}
\subseteq \cdots$ of $SCA$-subrings of $T$ such that $IT \cap R_{n}
\subseteq IR_{n+1}$ for every finitely generated ideal $I$ of $R_{n}$.

We now claim that $S = \bigcup_{i = 1}^{\infty} R_{i}$ is the
desired $SCA$-subring.  To see this, first note that $S$ is indeed an
$SCA$-subring, and that $R \subseteq S \subseteq T$.  Now set $I =
(y_{1},y_{2}, \ldots , y_{k})S$ and let $c \in IT \cap S$.  Then
there exists an $N \in \mathbb{N}$ such that $c, y_{1}, \ldots,
y_{k} \in R_{N}$.  Thus $c \in (y_{1},\ldots,y_{k})T \cap R_{N}
\subseteq (y_{1},\ldots,y_{k})R_{N+1} \subset IS$.  From this it
follows that $IT \cap S = I$, so the fourth condition of the Lemma
holds.
\end{proof}

We now construct a domain $A$ that has the desired completion, as
well as other interesting properties.

\begin{lem} \label{compto}
Let $(T,M)$ be a complete local ring of dimension at least one, and
$G$ a set of nonmaximal prime ideals of $T$ where $G$ contains the
associated primes of $T$ and such that the set of maximal elements of $G$,
call
it $C$, is
finite.  Moreover suppose that if $q \in \mbox{Spec}T$ with $q
\subseteq P$ for some $P \in G$
then $q \in G$. Also suppose that
for each prime ideal $P \in G$, $P$ contains no nonzero integers of
$T$.  Then
there exists a local domain $A$ such that
\begin{enumerate}
\item $\hat A = T$
\item If $p$ is a nonzero prime ideal of $A$, then $T \otimes_{A}
k(p) \cong k(p)$ where $k(p) = A_{p}/pA_{p}$
\item The generic formal fiber of $A$ is exactly $G$ (and so has
maximal ideals $C$).
\item If $I$ is a nonzero ideal of $A$, then $A/I$ is complete.
\end{enumerate}
\end{lem}

We note here that although the second and fourth conditions
of this lemma may not seem relevant, they will prove
useful later when, under certain circumstances, we show
that $A$ can be forced to be excellent.

\begin{proof}
The proof is quite similar to Lemma
8 in \cite{loex}.
Define
\[
\Omega = \{ u + J \in T/J | J \mbox{ is an ideal of } T \mbox{ with }
J \not\subseteq P \mbox{ for every } P \in G\}
\]
We claim that $|\Omega | \leq |T|$.  Since $T$ is infinite and
Noetherian, $|\{ J \mbox{ is an ideal of } T \mbox{ with } J
\not\subseteq P \mbox{ for all } P \in G\} | \leq |T|$.  Now, if $J$ is an
ideal of $T$, then $|T/J| \leq |T|$.  It follows that $|\Omega | \leq
|T|$.

Well order $\Omega$ so that each element has fewer than $|\Omega |$
predecessors.  Let $0$ denote the first element of $\Omega$.  Define
$R'_{0}$ to be the prime subring of $T$, and let $R_{0}$ simply
denote $R_{0}'$ localized at $R'_{0} \cap M$.  Note that $R_{0}$ is an
$SCA$-subring.

Now recursively define a family of $SCA$-subrings as follows, starting
with $R_{0}$.  Let $\lambda \in \Omega$ and
assume that $R_{\beta}$ has already been defined for all $\beta <
\lambda$.  Then $\gamma(\lambda) = u + J$ for some ideal $J$ of $T$
with $J \not\subseteq P$ for all $P \in G$ and thus all $P \in C$.
If $\gamma(\lambda) <
\lambda$, use Lemma \ref{sca} to obtain an $SCA$-subring $R_{\lambda}$
such that $R_{\gamma(\lambda)} \subseteq R_{\lambda} \subseteq T$,
$u + J \in \mbox{Image}(R_{\lambda} \to T/J)$ and for every finitely
generated ideal $I$ of $R_{\lambda}$ the property $IT \cap
R_{\lambda} = I$ holds.  Moreover, this gives us that $R_{\lambda}
\cap J \neq (0)$.  If $\gamma(\lambda) = \lambda$,
define $R_{\lambda} = \bigcup_{\beta < \lambda}R_{\beta}$.  Then we
have $R_{\lambda}$ is an $SCA$-subring for all $\lambda \in \Omega$.
We claim that $A = \bigcup_{\lambda \in \Omega}R_{\lambda}$ is the
desired domain.

We will first show that the generic formal fiber ring of $A$ has the
desired properties.  As each $R_{\lambda}$ is an $SCA$-subring, we
have $R_{\lambda} \cap P = (0)$ for each $P \in C$ and thus each $P
\in G$.  Therefore $A \cap P =
(0)$ for each $P \in G$ as well.  Moreover, if $J$ is an ideal of $T$
with $J \not\subseteq P$ for all $P \in G$, then $0 + J \in \Omega$.
Therefore, $\gamma(\lambda) = 0 + J$ for some $\lambda \in \Omega$
with $\gamma(\lambda) < \lambda$.  By construction, $R_{\lambda} \cap
J \neq (0)$.  It follows that $J \cap A \neq (0)$.  Hence the generic
formal fiber of $A$ is exactly $G$, and has maximal ideals $C$.

Now we show that the completion of $A$ is $T$.  To do this, we will
use Proposition \ref{prop1}.  Note that as each prime ideal $P \in G$
is nonmaximal in $T$, we have that $M^{2}$ is not contained in any $P
\in G$.  Thus by the construction, the map $A \to T/M^{2}$ is
surjective.  Now let $I$ be a finitely generated ideal of $A$ with
$I = (y_{1}, \ldots, y_{k})$.  Let $c \in IT \cap A$.  Then $\{ c,
y_{1}, \ldots, y_{k}\} \subseteq R_{\lambda}$ for some $\lambda \in
\Omega$ with $\gamma(\lambda) < \lambda$.  Again by the
construction, $(y_{1}, \ldots, y_{k})T \cap R_{\lambda} = (y_{1},
\ldots , y_{k})R_{\lambda}$.  As $c \in (y_{1},\ldots,y_{k})T \cap
R_{\lambda}$, we have that $c \in (y_{1},
\ldots,y_{k})R_{\lambda}\subseteq I$.  Hence $IT \cap A = I$ as
desired, and it follows that $A$ is Noetherian and its
completion is $T$.

To show the fourth condition is fairly simple.  Suppose that $I$ is a
nonzero ideal of $A$, and let $J = IT$.  If $J \subseteq P$ for some
$P \in G$, then $I \subseteq J \cap A \subseteq P \cap A = (0)$, a
contradiction.  Thus $J \not\subseteq P$ for every $P \in G$.  It
follows by construction that the map $A \to T/J$ is surjective.
Now since $A \cap J = A \cap IT = I$,
the map $A/I \to T/J$ is an isomorphism, making $A/I$ complete.

Finally, we prove the second condition.
Let $p$ be a nonzero prime ideal of $A$.  Then $A/p$ is
complete, so we have $T \otimes_{A}
k(p) \cong (T/pT)_{\overline{A-p}} \cong (A/p)_{\overline{A-p}} \cong
A_{p}/pA_{p} = k(p)$, as desired.
\end{proof}

\section{\textbf{\sc{The Main Theorem and Corollaries}}}
\label{nonexcellent}


For the next proof, we will in fact only need two of the previous
Lemma's four results,
namely the first one and the third one.  We are finally ready to
arrive at our main theorem.

\begin{thm} \label{gffthm}
Let $(T,M)$ be a complete local ring, and $G \subseteq \mbox{Spec}T$
such that $G$ is nonempty and
the number of maximal elements of $G$ is finite.
Then there exists a local domain $A$ such
that $\hat A = T$ and the generic formal fiber of $A$ is exactly $G$
if and only if $T$ is a field and $G = \{(0)\}$
or the following conditions hold.
\begin{enumerate}
\item $M \notin G$, and $G$ contains all the associated primes of $T$
\item If $Q \in G$ and $P\in \mbox{Spec}T$ with $P \subseteq Q$ then $P \in G$
\item If $Q \in G$ then $Q \cap \mbox{ prime subring of }T = (0)$
\end{enumerate}
\end{thm}

\begin{proof}
First, the forward direction.
Suppose that $T$ is not a field, and that there exists a local
domain $A$ such that $\hat A = T$ and the generic formal fiber of $A$
is exactly $G$.  Suppose that $M \in G$.  But then from our
assumptions $M \cap A = (0)$, which implies that the maximal ideal of
$A$ is $(0)$ and thus the maximal ideal of $T$ is zero as well, which
implies that $T$ is a field, a contradiction.  Thus $M \notin G$.
Moreover, if $G$ does not contain all of the associated primes of $T$,
then, as the generic formal fiber of $A$ is exactly $G$, there
must be an associated prime ideal $P$ of $T$ such that $P
\cap A \neq (0)$, a contradiction.  Therefore $G$ contains all of the
associated primes of $T$.

That the second requirement holds true is clear, as if one ideal $Q$ is
in the generic formal fiber of $A$, then $Q \cap A = (0)$ and thus if
$P \subseteq Q$ then $P \cap A = (0)$ and $P$ is also in the generic
formal fiber of $A$ and thus is contained in $G$, as desired.

In order to see that the intersection of each $P \in G$ with the prime
subring of $T$ is $(0)$, note that $P \cap A = (0)$ for each $P \in
G$.  Since $T$ is the completion of an integral domain $A$, the
unity element of $T$ must be in $A$.  Hence the prime subring of
$T$ is also in $A$, and $P$ contains no nonzero integers of $T$ for
all $P \in G$.
%

On the other hand, suppose that $T$ is a field.  Then the only prime ideal
of $T$ is $(0)$, and consequently $G = \{ (0)\}$.  Thus, since as
a field $T$ is a completion of itself, $A = T$ and the generic formal
fiber of $A$ is $\{ (0)\} = G$ as desired, so we're done.

Now we prove the backwards direction.  If $T$ is a field, then
$A = T$ works.
So suppose that $T$ is not a field and that all the above conditions
hold.  The first condition gives us that dim$T \geq 1$.
Now, use Lemma \ref{compto} to construct
the desired domain $A$.
\end{proof}

\begin{ex}
\textit{Let $T = \mathbb{C}[[x,y,z]]$, $C = \{ (x,y),(z)\}$, and $G$ be those
prime ideals $P$ of $T$ such that $P \subseteq (x,y)$ or $P
\subseteq (z)$.  Is there a local
domain $A$ such that $\hat A = T$ and the generic formal fiber of $A$
is exactly $G$ and has maximal ideals the elements of $C$?}

Clearly $T$ is local, with maximal ideal $M = (x,y,z)$, and $C$ is
finite.  Thus we may use
Theorem \ref{gffthm}.  Certainly $M \notin G$.  Moreover, $T$ is a
domain, so it has no zero divisors, and hence no associated primes
other than $(0)$, which is in $G$.
The way in which we defined $G$ makes it evident that the second
condition of the Theorem holds.  The third condition is also
easy to see since the prime subring of $T$ is
$\mathbb{Z}$ and all integers are
units.
Thus $T$ satisfies the three conditions of the Corollary, and hence
there exists a domain $A$ such that $\hat A = \mathbb{C}[[x,y,z]]$
and the generic formal fiber of $A$ is exactly $G$ with maximal
ideals $(x,y)$ and $(z)$ as desired.
\end{ex}

We now state the local version of Theorem \ref{gffthm} in
the following corollary.

\begin{cor} \label{localthm}
Let $(T,M)$ be a complete local ring and $P$ a prime ideal of $T$.  Then there
exists a local integral domain $A$ such that $\hat{A} = T$ and the
generic formal fiber of $A$ is local with maximal ideal $P$ if and only if
either $T$ is a field and $P = (0)$ or the following two conditions hold:
\begin{enumerate}
  \item $P$ is nonmaximal in $T$ and contains all the associated prime ideals
    of $T$
  \item $P \cap$ prime subring of $T = (0)$

\end{enumerate}

\end{cor}

\begin{proof}
That the forward direction holds true is obvious from Theorem
\ref{gffthm}, letting $C = \{ P\}$.

To see the backwards direction in the case where $T$ is not a field
(if $T$ is a field, then $A = T$ works), $P$ is a nonmaximal
ideal of $T$ containing all the
associated primes of $T$, and $P \cap \mbox{ prime subring of } T =
(0)$, simply let the set $G = \{ Q \in \mbox{Spec}T \ : \ Q \subseteq
P\}$.  It follows from Theorem \ref{gffthm} that there exists a domain
$A$ such that $\hat A = T$ and the generic formal fiber of $A$ is
local with maximal ideal $P$.
\end{proof}

\begin{ex}
\textit{Let $T = \frac{\mathbb{C}[[x,y,z]]}{(xy)}$.
Does there exist a local integral domain $A$ such that $\hat
A = T$ and the generic formal fiber of $A$ is local with maximal
ideal $(x,y)$?}

Let $P = (x,y)$.  It is easy to see that $P$ satisfies the
two conditions of Corollary \ref{localthm}.  So,
there exists a local domain $A$ such that $\hat A =
\frac{\mathbb{C}[[x,y,z]]}{(xy)}$ and
and the generic formal fiber of $A$ is local with maximal ideal
$(x,y)$.
\end{ex}

\section{\textbf{\sc{The Excellent Case}}} \label{excellent}


We now consider under what conditions the ring $A$ can be
made excellent.  Using
the same building blocks as the previous theorem, we are able to come up
with a characterization in the characteristic
zero case of those complete local rings that are the
completion of a local excellent domain possessing a specific
generic formal fiber.

\begin{thm} \label{exthm}
Let $(T,M)$ be a complete local ring containing the integers.  Let $G
\subseteq \mbox{Spec}T$ such that $G$ is nonempty and the number of
maximal elements of $G$ is finite.  Then there exists an
excellent local domain $A$ with $\hat{A} = T$ and
such that $A$ has generic formal fiber exactly $G$ if
and only if $T$ is a field and $G = \{(0)\}$
or the following conditions hold.
\begin{enumerate}
\item $M \notin G$, and $G$ contains all the associated primes of $T$
\item If $Q \in G$ and $P\in \mbox{Spec}T$ with $P \subseteq Q$ then $P \in G$
\item If $Q \in G$, then $Q \cap \mbox{ prime subring of }T = (0)$
\item $T$ is equidimensional
\item $T_{P}$ is a regular local ring for all maximal elements $P \in
G$.
\end{enumerate}
\end{thm}

\begin{proof}
Assume that $T$ is the completion of an excellent domain $A$ having
generic formal fiber exactly $G$ with maximal ideals the maximal
elements of $G$.  If $\mbox{dim}T = 0$ then $T$ is a field
and $G = \{(0)\}$.  Thus consider the case
where $\mbox{dim}T \geq 1$.

As $A$ is excellent, it is universally catenary.
Hence, $A$ is formally catenary and it follows that $A/(0) \cong
A$ is formally equidimensional.  Thus the completion, $T$, is
equidimensional.

%
%
The first three conditions can be shown by the exact same arguments as
the three conditions in
Theorem \ref{gffthm}.

To see that the fifth condition holds, note that
 the maximal ideals of $T \otimes_{A} k(0)$ are the
maximal elements of $G$.  Let $P$ be one of these maximal elements.
Then $T \otimes_{A} k(0)$ localized at $P$ is isomorphic to $T_{P}$.
Since $A$ is excellent, $T
\otimes_{A} k(0)$ is regular, implying that $T_{P}$ is a regular local
ring for every maximal element of $G$ as desired.

Conversely, first suppose that $T$ is a field and
$G = \{(0)\}$.  Then $A = T$ works.
So, suppose that $T$ is not a field and
that all of the five conditions hold true for
some complete local ring $T$ and some nonempty set $G$ of prime ideals
of $T$ such that the number of maximal elements of $G$ is finite.
We want to show that there exists an excellent domain $A$ possessing
generic formal fiber exactly $G$.  Note that conditions $(1)$ and
$(5)$ imply that $T$ is reduced by the following argument.  Suppose
that these two conditions are true, but that $T$ is not reduced.
Then there exists a nonzero $x \in T$ such that $x^{n} = 0$ in $T$.  Now
consider the ideal $(0 : x)$ of $T$.  Now, $(0 :
x) \subseteq Q_{1} \cup Q_{2} \cup \cdots \cup Q_{n}$, where the
$Q_{i}$ are the associated prime ideals of $T$.  But then by the Prime
Avoidance Theorem $(0 : x) \subseteq Q_{i}$ for some $i$.  Moreover,
as $G$ contains the associated primes of $T$, $Q_{i} \subseteq P$ for
some $P$ a maximal element of $G$.  Consider the regular local ring
$T_{P}$, and note that regular local rings are domains.
Now $\frac{x^{n}}{1} = \frac{0}{1}$ in $T_{P}$.  Thus, as
$T_{P}$ is a domain, $\frac{x}{1} = \frac{0}{1}$ in $T_{P}$.  But then
there exists an element $s \notin P$ such that $sx = 0$.  This implies
that $s \in (0 : x)$, however, which indicates that $s \in P$, a
contradiction.  Therefore, $T$ must be reduced as desired.

Now if $\mbox{dim}T = 0$ then
 $T$ is a field and we're in the first case.  Suppose,
on the other hand, that $\mbox{dim}T \geq 1$.  Then use Lemma
\ref{compto} to construct the domain $A$.  We claim that $A$ is
excellent with generic formal fiber exactly $G$.  From the
construction of $A$, $A$ has the desired generic formal fiber.  To see
that $A$ is excellent, suppose that $p$ is a nonzero prime ideal of
$A$.  Then from Lemma \ref{compto} we have $T \otimes_{A} k(p) \cong k(p)$.
Now let $L$ be a finite field extension of $k(p)$.  Then $T \otimes_{A}
L \cong T \otimes_{A} k(p) \otimes_{k(p)} L \cong k(p) \otimes_{k(p)}
L \cong L$.  Thus the fiber over $p$ is geometrically regular.  Now
$T_{P}$ is regular by assumption for every maximal element $P$ of
$G$.  It follows that $T \otimes_{A} k(0)$ is regular.  Now since $T$
contains the integers, so does $A$.  It follows that $k((0))$ is a
field of characteristic zero, and hence that $T \otimes_{A} L$ is
regular for every finite field extension $L$ of $k((0))$.  Thus all
of the formal fibers of $A$ are geometrically regular.  Since $A$ is
formally equidimensional it is universally catenary,
and thus $A$ is excellent.  Hence $A$ is the
desired domain.
\end{proof}

Notice that this proof fails if the
characteristic of $T$ is $p > 0$, as the $A$ we construct may not have a
geometrically regular generic formal fiber.  It is worth
noting, though, that all the
other fibers are geometrically regular and so the only
obstruction to $A$ being excellent is that the generic formal
fiber may not be geometrically regular.

We now state the local version of Theorem \ref{exthm}.
 Arguably more elegant than the
previous theorem, this more specific theorem has fewer conditions, and
is thus may prove to be more practical.

\begin{cor} \label{exlocal}
Let $(T,M)$ be a complete local ring containing the integers and $P$ a
prime ideal of $T$.
Then $T$ is the completion of a local excellent domain $A$ possessing a
local
generic formal fiber with maximal ideal $P$ if and
only if $T$ is a field and $P = (0)$
or the following three conditions hold:
\begin{enumerate}
	\item $P$ contains the associated prime ideals of $T$
	\item $P$ is a nonmaximal prime ideal of $T$ such that
	$P$ contains no nonzero integers of $T$
	\item $T_{P}$ is a regular local ring
\end{enumerate}
\end{cor}

\begin{proof}
The forward direction of this proof follows immediately from Theorem
\ref{exthm}.

Conversely, if $T$ is a field and $P = (0)$, then $A = T$ works.

So suppose $P$ is a nonmaximal prime ideal
of $T$ containing all the associated primes of $T$ such that $P \cap
\mbox{ prime subring of }T = (0)$, and
$T_{P}$ is a regular local ring.  It interesting to note that these
conditions alone imply that $T$ is an integral domain by the following
argument.
Suppose that $xy = 0$ in $T$ for some $x,y \in T$, $x,y \neq 0$.  But then
both $x$ and $y$ are zerodivisors of $T$, so $x,y \in P$.
 But now since $T_{P}$ is a domain,
then either $x/1 = 0/1$ or $y/1 = 0/1$.  WLOG, assume that $x/1 = 0/1$.
This implies that there
is some $s \not\in P$ such that $sx = 0$ in $T$.  But then $s$ is a zero
divisor or $T$, and hence $s \in P$, a contradiction.  Therefore, $T$
contains no nonzero zero divisors and $T$ is an integral domain.
Observe that $T$ being an integral domain
implies that $T$ is both reduced and equidimensional.

Now if
dim$T = 0$, then since $T$ is reduced $M$ must be an associated prime of $T$,
and moreover it is the only associated prime ideal of $T$.
Hence $T$ is a field and we're in the first case.

%
If $\mbox{dim}T \geq 1$,
let $G = \{ Q \in
\mbox{Spec}T \ : \ Q \subseteq P\}$, and
use Lemma \ref{compto} to construct
the domain $A$.  It is trivial to verify that the five conditions of
Theorem \ref{exthm} hold.
Theorem \ref{exthm} tells us that $A$ is then
excellent with generic formal fiber exactly $G$.  But the maximal
ideal of $G$ is $P$ by definition, so $A$ is excellent with local
generic formal with maximal ideal $P$ as desired.
%
\end{proof}

\begin{ex} \label{exeg}
\textit{Consider the complete local ring $T =
\frac{\mathbb{C}[[x,y,z]]}{(xy)}$ and $G = \{ Q \in
\mbox{Spec}T \ : \ Q \subseteq (x) \mbox{ or } Q \subseteq (y,z)\}$
a set of prime ideals of $T$ with maximal elements $(x)$ and $(y,z)$.
It is not difficult to check that $T$ and $G$ satisfy
the conditions of Theorem \ref{exthm} and so $T$ is
the completion of an excellent
local domain $A$ with generic formal fiber
exactly $G$.}

Note here that $T$ in the
above example is not a domain, and thus we should
not be able to find a prime ideal $P$ of
$T$ such that there exists an excellent domain $A$ that completes to
$T$ with local generic formal fiber with maximal ideal $P$.
(Recall that in the proof of Corollary \ref{exlocal}, we showed that
such a $T$ is necessarily an integral domain.)  Indeed
we can observe that this is true by seeing that any ideal we might choose
will either not contain $\mbox{Ass}T = \{(x),(y)\}$, or if it does
then $T_{P}$ will not be a regular local ring.
\end{ex}





\bigskip

\begin{center}
\begin{tabular}{cc}
P. Charters  &  S. Loepp \\
Department of Mathematics  &   Department of Mathematics and Statistics  \\
University of Texas at Austin  &  Williams College\\
Austin, TX  78712  &  Williamstown, MA  01267 \\
pcharter@wso.williams.edu  &  sloepp@williams.edu    \\
\end{tabular}
\end{center}

\end{document}